\documentclass[12pt,leqno]{amsart}
\usepackage{amsmath}
\usepackage{amsthm}
\usepackage{amssymb}
\def\overset#1#2{{\mathrel{\mathop {{#2}_{}}\limits^{#1}}}}
\def\underset#1#2{{\mathrel{\mathop {{}_{} {#2}}\limits_{{#1}_{}}}}}
\def\upplim_#1{\underset{#1}{\overline\lim}\;}
\def\lowlim_#1{\underset{#1}{\underline\lim}\;}
\newtheorem{lemma}[equation]{Lemma}

\newtheorem{theorem}[equation]{Theorem}

\newcommand{\C}{{\mathbf{C}}}

\newcommand{\N}{\mathbf{N}}

\renewcommand{\P}{{\mathbf{P}}}

\newcommand{\R}{{\mathbf{R}}}

\newcommand{\zero}{\mathrm{Zero}}

\newcommand{\Z}{\mathbf{Z}}

\numberwithin{equation}{section}

\begin{document}
\title[Second main theorems and uniqueness problem of meromorphic mappings]{Second main theorems and uniqueness problem of meromorphic mappings with moving hypersurfaces} 
\date{} 
\author{Si Duc Quang}

\maketitle

\begin{abstract} In this article, we establish some new second main theorems for meromorphic mappings of $\C^m$ into $\P^{n}(\C)$ and moving hypersurfaces with truncated counting functions. A uniqueness theorem for these mappings sharing few moving hypersurfaces without counting multiplicity is also given. This result is an improvement of the recent result of Dethloff - Tan \cite{DT02}. Moreover the meromorphic mappings in our result may be algebraically degenerate. The last purpose of this article is to study uniqueness problem in the case where the meromorphic mappings agree on small identical sets.
\end{abstract}

\def\thefootnote{\empty}

\setlength{\baselineskip}{16pt}
\maketitle

\section{Introduction}

In 2004, Min Ru \cite{MR} showed a second main theorem for algebraically nondegenerate meromorphic mappings and a family of hypersurfaces in weakly general position. After that, with the same assumptions, T. T. H. An and H. T. Phuong \cite{AP} improved the result of Min Ru by giving an explicit truncation level for counting functions. 

Recently, in  \cite{DT01} Dethloff and Tan generalized and improved the second main theorems of Min Ru and An - Phuong to the case of moving hypersurfaces. They proved that

\vskip0.2cm
\noindent
\textbf{Theorem A} (Dethloff - Tan \cite{DT01})
{\it Let $f$ be a nonconstant meromorphic map of $\C^m$ into $\P^n(\C)$. Let $\{Q_i\}_{i=1}^q$ be a set of slow (with respect to $f$) moving hypersurfaces in weakly general position with $\deg Q_j = d_j\ (1\le i\le q).$ Assume that $f$ is algebraically nondegenerate over $\tilde{\mathcal K}_{\{Q_i\}_{i=1}^q}$.  Then for any $\epsilon >0$ there exist positive integers $L_j\ (j = 1,....,q)$, depending only on $n,\epsilon$ and $d_j\ (j = 1,...,q) $ in an explicit way such that 
$$  ||\ (q-n-1-\epsilon)T_f(r)\le \sum_{i=1}^{q}\dfrac{1}{d_i}N^{[L_j]}_{Q_i(f)}(r)+o(T_f(r)) .$$}
\vskip0.2cm
Here, the truncation level $L_j$ is estimated by
$$ L_j\le \dfrac{d_j\cdot \binom{n+M}{n}t_{p_0+1}-d_j}{d}+1, $$
where $d$ is the least common multiple of  the $d_j'$s, $d=lcm (d_1,...,d_q)$, and
\begin{align*}
M=&d\cdot [2(n+1)(2^n-1)(nd+1)\epsilon^{-1}+n+1],\\ 
p_0=&[\dfrac{({\binom{n+M}{n}}^2\cdot\binom{q}{n}-1)\cdot\log ({\binom{n+M}{n}}^2\cdot\binom{q}{n}) }{\log (1+\dfrac{\epsilon}{2\binom{n+M}{n}M})}+1] ^2,\\
and  \ t_{p_0+1}<& \left ({\binom{n+M}{n}}^2\cdot\binom{q}{n}+p_0\right )^{\bigl ({\binom{n+M}{n}}^2\cdot\binom{q}{n}-1\bigl )},
\end{align*}
where $[x] = \max\{k\in \Z\ ;\ k\le x\}$ for a real number $x$.

By using this second main theorem, Dethloff and Tan proved a uniqueness theorem for meromorphic mappings which share slow moving hypersurfaces as follows.

Let $f,g:\C^m\rightarrow \P^n(\C)$ be two meromorphic mappings. Let $\{Q_i\}_{i=1}^q$ be $q$ moving hypersurfaces of $\P^n(\C)$ in weakly general position, $\deg Q_i=d_i,$ and let $d,d^*,\tilde d$ be respectively the least common multiple, the maximum number and the minimum number of the ${d_j}'$s. Take $M, p_0$ be as above with $\epsilon =1$ and set
\begin{align*}
t_{p_0+1}=& \left ({\binom{n+M}{n}}^2\cdot\binom{q}{n}+p_0\right )^{\bigl ({\binom{n+M}{n}}^2\cdot\binom{q}{n}-1\bigl )},\\
L=&[\dfrac{d^*\cdot \binom{n+M}{n}t_{p_0+1}-d^*}{d}+1].
\end{align*}
With the above notations, in 2011, Dethloff and Tan proved the following.

\vskip0.2cm
\noindent
\textbf{Theorem B} (Theorem 3.1 \cite{DT02}).\
{\it a) Assume that $f$ and $g$ are algebraically nondegenerate over $\tilde {\mathcal K}_{\{Q_j\}}$ such that:
$$i) \mathcal D^\alpha (\dfrac{f_k}{f_s})=\mathcal D^\alpha (\dfrac{g_k}{g_s})\text{ on } \bigcup_{i=1}^q(\zero Q_i(f)\cup \zero Q_i(g)),$$
for all $|\alpha|<p$, $p\in\Z^+$ and $0\le k\ne s\le n$.
 
Then for $q>n+\frac{2nL}{p\tilde d}+\frac{3}{2},$ we have $f\equiv g.$

b) Assume $f$ and $g$ as a) satisfy i) and
$$ \dim \bigl (\bigcap_{j=0}^{n}\zero Q_{i_j}(f)\bigl )\le m-2\ \forall 1\le i_0<\cdots <i_n\le q. $$
Then for $q>n+\frac{2L}{p\tilde d}+\frac{3}{2},$ we have $f\equiv g.$}

\noindent
\vskip0.2cm
However,  the number of moving hypersurfaces in Theorem B is still big, since the truncation levels given in Theorem A is far from the sharp. 

We also would like to note that, in all mentioned results on second main theorem of Min Ru, An - Phuong and Dethloff - Tan the algebraically nondegeneracy condition of the meromorphic mappings can not be removed and it plays an essential role in their proofs.

The first purpose of the present paper is to show some new second main theorems for meromorphic mappings sharing slow moving hypersurfaces with better truncation levels for counting functions. Moreover the mappings may be algebraically degenerate. Namely, we prove the following theorems.

\begin{theorem}\label{1.1} 
Let $f$ be a meromorphic mapping of $\C^m$ into $\P^n(\C)$. Let $Q_i\ (i=1,...,q)$ be slow (with respect to $f$) moving hypersurfaces of $\P^n(\C)$ in weakly general position with $\deg Q_i=d_i$, $q\ge nN+n+1,$ where $N=\binom{n+d}{n}-1$ and $d=lcm (d_1,...,d_q)$. Assume that $Q_i(f)\not\equiv 0\ (1\le i\le q)$. Then we have
$$ ||\ \dfrac{q}{nN+n+1}T_f(r)\le \sum_{i=1}^{q}\dfrac{1}{d_i}N^{[N]}_{Q_i(f)}(r)+o(T_f(r)).$$
\end{theorem}

\begin{theorem}\label{1.2} 
Let $f$ be a meromorphic mapping of $\C^m$ into $\P^n(\C)$. Let $Q_i\ (i=1,...,q)$ be slow (with respect to $f$) moving hypersurfaces of $\P^n(\C)$ in weakly general position with $\deg Q_i=d_i$, $q\ge N+2,$ where $N=\binom{n+d}{n}-1$ and $d=lcm (d_1,...,d_q)$. Assume that $f$ is algebraically nondegenerate over  $\mathcal {\tilde K}_{{\{Q_i\}}_{i=1}^q}$. Then we have
$$ ||\ \dfrac{q}{N+2}T_f(r)\le \sum_{i=1}^{q}\dfrac{1}{d_i}N^{[N]}_{Q_i(f)}(r)+o(T_f(r)).$$
\end{theorem}

The second purpose of this paper is to show a uniqueness theorem for meromorphic mappings sharing slow moving hypersurfaces without counting multiplicity. We will prove the following.

\begin{theorem}\label{1.3}
Let $f$ and $g$ be nonconstant meromorphic mappings of $\C^m$ into $\P^n(\C)$. Let $Q_i\ (i=1,...,q)$ be a set of slow (with respect to $f$ and $g$) moving hypersurfaces in $\P^n(\C)$ in weakly general position with $\deg Q_i=d_i$. Put $d=lcm (d_1,...,d_{n+2})$ and $N=\binom{n+d}{n}-1$. Let $k\ (1\le k\le n)$ be an integer. Assume that 

(i) $\dim \bigl (\bigcap_{j=0}^{k}\zero Q_{i_j}(f)\bigl )\le m-2$ for every $1\le i_0<\cdots <i_k\le q,$

(ii) $f=g$ on $\bigcup_{i=1}^q\bigl (\zero Q_i(f)\cup \zero Q_{i_j}(g)\bigl )$.\\
Then the following assertions hold:

a) If $q>\dfrac{2kN(nN+n+1)}{d}$ then $f=g.$

b) In addition to the assumptions (i)-(ii), we assume further that both $f$ and $g$ are algebraically nondegenerate over 
$\mathcal{\tilde K}_{\{Q_i\}_{i=1}^q}$. If $q>\dfrac{2kN(N+2)}{d}$, then $f=g.$
\end{theorem}
We note that the numbers of hypersurfaces in our results are really reduced when compared to that in Theorem B of Dethloff - Tan. Also by introducing some new techniques, we simplify their proofs.

We would like to emphasize  here that in all Theorem \ref{1.3} and previous results on the uniqueness problem, the meromorphic mappings always are assumed to agree on the \textit{"inverse images"} of all moving hypersurfaces. Our last purpose in this paper is to show an algebraic relation between meromorphic mappings in the case where they agree on the \textit{"inverse images"} of only $n+2$ moving hypersurfaces. Namely, we will prove the following. 

\begin{theorem}\label{1.4}
Let $f$ and $g$ be nonconstant meromorphic mappings of $\C^m$ into $\P^n(\C)$. Let $Q_i\ (i=1,...,q)$ be a set of slow (with respect to $f$ and $g$) moving hypersurfaces in $\P^n(\C)$ in weakly general position with $\deg Q_i=d_i$. Put $d=lcm (d_1,...,d_q)$,\\
$ L_j=\bigl [ \dfrac{d_j\cdot \binom{n+M}{n}t_{p_0+1}-d_j}{d}+1\bigl ]$, where 
$M=d\cdot \left (4(n+1)(2^n-1)(nd+1)+n+1\right ), p_0= \bigl [\dfrac{({\binom{n+M}{n}}^2\cdot\binom{q}{n}-1)\cdot\log ({\binom{n+M}{n}}^2\cdot\binom{q}{n}) }{\log (1+\dfrac{1}{4\binom{n+M}{n}M})}+1\bigl ] ^2$
and $t_{p_0+1}=\left ({\binom{n+M}{n}}^2\cdot\binom{q}{n}+p_0\right )^{\bigl ({\binom{n+M}{n}}^2\cdot\binom{q}{n}-1\bigl )}.$
 Assume that $f$ and $g$ are algebraically nondegenerate over $\mathcal{\tilde K}_{\{Q_i\}_{i=1}^q}$ and 

(i) $\dim \bigl (\bigcap_{j=0}^{k}\zero Q_{i_j}(f)\bigl )\le m-2$ for every $1\le i_0<\cdots <i_k\le n+2,$

(ii) $\min\{\nu^0_{Q_i(f)}(z), L_i\}=\min\{\nu^0_{Q_i(g)}(z), L_i\}$ for every $n+3\le i\le q,$

(ii) $f=g$ on $\bigcup_{i=1}^{n+2}\bigl (\zero Q_i(f)\cup \zero Q_{i_j}(g)\bigl )$.

\noindent
If  $q\ge n+2+2kL$, where $L=\max_{1\le i\le n+2}\frac{dL_i}{d_i}$ then there exist at least $[\frac{q-n-2}{2}]+1$ indices $n+3\le i_1<\cdots <i_{[\frac{q-n-2}{2}]+1}$ such that
$$ \dfrac{Q_{i_1}(f)}{Q_{i_1}(g)}=\dfrac{Q_{i_2}(f)}{Q_{i_2}(g)}=\cdots =\dfrac{Q_{i_{[\frac{q-n-2}{2}]+1}}(f)}{Q_{i_{[\frac{q-n-2}{2}]+1}}(g)}.$$
\end{theorem}

{\bf Acknowledgements.} This work was supported in part by a NAFOSTED grant of Vietnam.

\section{Basic notions and auxiliary results from Nevanlinna theory}

\noindent
{\bf 2.1.}\ We set $||z|| = \big(|z_1|^2 + \dots + |z_m|^2\big)^{1/2}$ for
$z = (z_1,\dots,z_m) \in \C^m$ and define
\begin{align*}
B(r) := \{ z \in \C^m : ||z|| < r\},\quad
S(r) := \{ z \in \C^m : ||z|| = r\}\ (0<r<\infty).
\end{align*}
Define 
$$v_{m-1}(z) := \big(dd^c ||z||^2\big)^{m-1}\quad \quad \text{and}$$
$$\sigma_m(z):= d^c \text{log}||z||^2 \land \big(dd^c \text{log}||z||^2\big)^{m-1}
 \text{on} \quad \C^m \setminus \{0\}.$$
{\bf 2.2.}\ Let $F$ be a nonzero holomorphic function on a domain $\Omega$ in $\C^m$. For a set $\alpha = (\alpha_1,...,\alpha_m) $ of nonnegative integers, we set $|\alpha|=\alpha_1+...+\alpha_m$ and 
$\mathcal {D}^\alpha F=\dfrac {\partial ^{|\alpha|} F}{\partial ^{\alpha_1}z_1...\partial ^{\alpha_m}z_m}.$
We define the map $\nu_F : \Omega \to \Z$ by
$$\nu_F(z):=\max\ \{k: \mathcal {D}^\alpha F(z)=0 \text { for all } \alpha \text { with }|\alpha|<k\}\ (z\in \Omega).$$

We mean by a divisor on a domain $\Omega$ in $\C^m$ a map $\nu : \Omega \to \Z$ such that, for each $a\in \Omega$, there are nonzero holomorphic functions $F$ and $G$ on a connected neighborhood $U\subset \Omega$ of $a$ such that $\nu (z)= \nu_F(z)-\nu_G(z)$ for each $z\in U$ outside an analytic set of dimension $\le m-2$. Two divisors are regarded as the same if they are identical outside an analytic set of dimension $\le m-2$. For a divisor $\nu$ on $\Omega$ we set $|\nu| := \overline {\{z:\nu(z)\ne 0\}},$ which
is either a purely $(m-1)$-dimensional analytic subset of $\Omega$ or an empty set.

Take a nonzero meromorphic function $\varphi$ on a domain $\Omega$ in $\C^m$. For each $a\in \Omega$, we choose nonzero holomorphic functions $F$ and $G$ on a neighborhood $U\subset \Omega$ such that $\varphi = \dfrac {F}{G}$ on $U$ and $\dim (F^{-1}(0)\cap G^{-1}(0))\le m-2,$ and  we define the divisors $\nu^0_\varphi,\ \nu^\infty_\varphi$ by 
$ \nu^0_\varphi := \nu_F,\  \nu^\infty_\varphi :=\nu_G$, 
which are independent of choices of $F$ and $G$ and so globally well-defined on $\Omega$.

\noindent
{\bf 2.3.}\ For a divisor $\nu$ on $\C^m$ and for a positive integer $M$ or $M= \infty$, we define the counting function of $\nu$ by
$$\nu^{[M]}(z)=\min\ \{M,\nu(z)\},$$
\begin{align*}
n(t) =
\begin{cases}
\int\limits_{|\nu|\,\cap B(t)}
\nu(z) v_{m-1} & \text  { if } m \geq 2,\\
\sum\limits_{|z|\leq t} \nu (z) & \text { if }  m=1. 
\end{cases}
\end{align*}

Similarly, we define \quad $n^{[M]}(t).$

Define
$$ N(r,\nu)=\int\limits_1^r \dfrac {n(t)}{t^{2m-1}}dt \quad (1<r<\infty).$$

Similarly, we define  \ $N(r,\nu^{[M]})$
and denote it by \ $N^{[M]}(r,\nu)$.

Let $\varphi : \C^m \longrightarrow \C $ be a meromorphic function.
Define
$$N_{\varphi}(r)=N(r,\nu^0_{\varphi}), \ N_{\varphi}^{[M]}(r)=N^{[M]}(r,\nu^0_{\varphi}).$$

For brevity we will omit the character $^{[M]}$ if $M=\infty$.

\noindent
{\bf 2.4.}\ Let $f : \C^m \longrightarrow \P^n(\C)$ be a meromorphic mapping.
For arbitrarily fixed homogeneous coordinates
$(w_0 : \dots : w_n)$ on $\P^n(\C)$, we take a reduced representation
$f = (f_0 : \dots : f_n)$, which means that each $f_i$ is a  
holomorphic function on $\C^m$ and 
$f(z) = \big(f_0(z) : \dots : f_n(z)\big)$ outside the analytic set
$\{ f_0 = \dots = f_n= 0\}$ of codimension $\geq 2$.
Set $\Vert f \Vert = \big(|f_0|^2 + \dots + |f_n|^2\big)^{1/2}$.

The characteristic function of $f$ is defined by 
\begin{align*}
T_f(r)= \int\limits_{S(r)} \log\Vert f \Vert \sigma_m -
\int\limits_{S(1)}\log\Vert f\Vert \sigma_m.
\end{align*}

\noindent
\textbf{2.5.}\ Let $\varphi$ be a nonzero meromorphic function on $\C^m$, which are occasionally regarded as a meromorphic map into $\P^1(\C)$. The proximity function of $\varphi$ is defined by
$$m(r,\varphi):=\int_{S(r)}\log \max\ (|\varphi|,1)\sigma_m.$$
The Nevanlinna's characteristic function of $\varphi$ is defined as follows
$$ T(r,\varphi):=N_{\frac{1}{\varphi}}(r)+m(r,\varphi). $$
Then 
$$T_\varphi (r)=T(r,\varphi)+O(1).$$
The function $\varphi$ is said to be small (with respect to $f$) if $||\ T_\varphi (r)=o(T_f(r))$.
Here, by the notation $''|| \ P''$  we mean the assertion $P$ holds for all $r \in [0,\infty)$ excluding a Borel subset $E$ of the interval $[0,\infty)$ with $\int_E dr<\infty$.

We denote by $\mathcal M$ (resp. $\mathcal K_f$) the field of all meromorphic functions (resp. small meromorphic functions) on $\C^m$.

\noindent
{\bf 2.6.}\ Denote by $\mathcal H_{\C^m}$ the ring of all holomorphic functions on $\C^m.$
Let $Q$ be a homogeneous polynomial in $\mathcal H_{\C^m}[x_0,\dots,x_n]$  of
degree $d \geq 1.$ Denote by $Q(z)$ the homogeneous  polynomial over $\C$ obtained by substituting a specific point $z \in \C^m$ into the coefficients of $Q$. We also call  a moving  hypersurface in $\P^n (\C )$  each homogeneous polynomial $Q \in\mathcal H_{\C^m}[x_0,\dots,x_n]$  such that the common zero set of all coefficients of $Q$ has codimension at least two.

Let $Q$ be a moving hypersurface in $\P^n(\C)$ of degree $d\ge 1$ given by
$$ Q(z)=\sum_{I\in\mathcal I_d}a_I\omega^I, $$
where $\mathcal I_d=\{(i_0,...,i_n)\in\N_0^{n+1}\ ;\ i_0+\cdots +i_n=d\}$, $a_I\in\mathcal H_{\C^m}$ and $\omega^I=\omega_0^{i_0}\cdots\omega_n^{i_n}$. We consider the meromorphic mapping $Q':\C^m\rightarrow\P^N(\C)$, where $N=\binom{n+d}{n}$, given by
$$ Q'(z)=(a_{I_0}(z):\cdots :a_{I_N}(z))\ (\mathcal I_d=\{I_0,...,I_N\}). $$
The moving hypersurfaces $Q$ is said to be "slow" (with respect to $f$) if $||\ T_{Q'}(r)=o(T_f(r))$. This is equivalent to $||T_{\frac{a_{I_i}}{a_{I_j}}}(r)=o(T_f(r))$ for every $a_{I_j}\not\equiv 0.$

Let $\{Q_i\}_{i=1}^q$ be a family of moving hypersurfaces in $\P^n(\C)$, $\deg Q_i=d_i$. Assume that
$$ Q_i=\sum_{I\in\mathcal I_{d_i}}a_{iI}\omega^I. $$
We denote by $\tilde{\mathcal K}_{\{Q_i\}_{i=1}^q}$ the smallest subfield of $\mathcal M$ which contains $\C$ and all $\frac{a_iI}{a_{iJ}}$ with $a_{iJ}\not\equiv 0$.  We say that $\{Q_i\}_{i=1}^q$ are in weakly general position if there exists $z \in \C^m$ such that all $a_{iI}\ (1\le i\le q,\ I\in\mathcal I)$ are holomorphic at $z$ and for any $1 \leq i_0 < \dots < i_n \leq q$ the system of equations
\begin{align*}
\left\{ \begin{matrix}
Q_{i_j}(z)(w_0,\dots,w_n) = 0\cr
0 \leq j \leq n\end{matrix}\right.
\end{align*}
has only the trivial solution $w = (0,\dots,0)$ in $\C^{n+1}$.

\noindent
\textbf{2.7.}\ Let $f$ be a nonconstant meromorphic map of $\C^m$ into $\P^n(\C)$. Denote by
$\mathcal C_f$ the set of all non-negative functions $h : \C^m\setminus A\longrightarrow [0,+\infty]\subset\overline\R$, which are of the form
$$ h=\dfrac{|g_1|+\cdots +|g_l|}{|g_{l+1}|+\cdots +|g_{l+k}|}, $$
where $k,l\in\N,\ g_1,...., g_{l+k}\in\mathcal K_f\setminus\{0\}$ and $A\subset\C^m$, which may depend on
$g_1,....,g_{l+k}$, is an analytic subset of codimension at least two. Then, for $h\in\mathcal C_f$ we have
$$\int\limits_{S(r)}\log h\sigma_m= o(T_f (r)).$$

\noindent
{\bf Lemma 2.8}\ (Lemma 2 \cite{DT01}).
{\it Let $\{Q_i\}_{i=0}^n$ be a set of homogeneous polynomials of degree $d$ in $\mathcal K_f [x_0,..., x_n]$. Then there exists a function $h_1\in\mathcal C_f$ such that, outside an analytic set of $\C^m$ of codimension at least two,
$$ \max_{i\in\{0,...,n\}}|Q_i(f_0,...,f_n)|\le h_1||f||^d .$$
If, moreover, this set of homogeneous polynomials is in weakly general position, then there
exists a nonzero function $h_2\in\mathcal C_f$ such that, outside an analytic set of $\C^m$ of
codimension at least two,
$$h_2||f||^d \le  \max_{i\in\{0,...,n\}}|Q_i(f_0,...,f_n)|.$$}
 
\noindent
{\bf 2.9. Lemma on logarithmic derivative} \ (Lemma 3.11 \cite{Shi}) .\ {\it Let $f$ be a nonzero meromorphic function on $\C^m.$ Then 
$$\biggl|\biggl|\quad m\biggl(r,\dfrac{\mathcal{D}^\alpha (f)}{f}\biggl)=O(\log^+T(r,f))\ (\alpha\in \Z^m_+).$$}

\noindent
{\bf 2.10.}\ Assume that $\mathcal {L}$ is a subset of a vector space $V$ over a field $\mathcal R$. We say that the set $\mathcal {L}$ is {\it  minimal} over $\mathcal R$ if it is linearly dependent over $\mathcal R$ and each proper subset of $\mathcal L$ is linearly independent over  $\mathcal {R}.$

Repeating the argument in (Prop. 4.5 \cite{Fu}), we have the following:

\noindent
{\bf Proposition 2.11.}\  {\it Let $\Phi_0,...,\Phi_k$ be meromorphic functions on $\C^m$ such that $\{\Phi_0,...,\Phi_k\}$ 
are  linearly independent over $\C.$ Then  there exists an admissible set  
$$\{\alpha_i=(\alpha_{i1},...,\alpha_{im})\}_{i=0}^k \subset \Z^m_+$$
with $|\alpha_i|=\sum_{j=1}^{m}|\alpha_{ij}|\le k \ (0\le i \le k)$ such that the following are satisfied:

(i)\  $\{{\mathcal D}^{\alpha_i}\Phi_0,...,{\mathcal D}^{\alpha_i}\Phi_k\}_{i=0}^{k}$ is linearly independent over $\mathcal M,$\ i.e., \ 
$\det{({\mathcal D}^{\alpha_i}\Phi_j)}\not\equiv 0.$ 

(ii) $\det \bigl({\mathcal D}^{\alpha_i}(h\Phi_j)\bigl)=h^{k+1}\cdot \det \bigl({\mathcal D}^{\alpha_i}\Phi_j\bigl)$ for
any nonzero meromorphic function $h$ on $\C^m.$}

\section{Second main theorems for moving hypersurfaces}

In order to prove Theorem \ref{1.1} we need the following.

\begin{lemma}\label{3.1} Let $f$ be as in Theorem \ref{1.1}. Let $\{Q_i\}_{i=0}^{n(N+1)}$ be a set of  homogeneous polynomials in $\mathcal K_f[x_0,...,x_n]$ of common degree $d$ in weakly general position, where $N=\binom{n+d}{n}-1$. Assume that $Q_i(f)\not\equiv 0\ (0\le i\le n(N+1))$. Then there exist a subset $B$ of  $\{Q_i(f)\ ;\ 0\le i\le n(N+1)\}$ and subsets  $I_1,...,I_k$ of  $B$ such that 
the following are satisfied:

$\mathrm(i)$\ \ $I_1$ is minimal, $I_i$ is independent over $\mathcal {K}_f\ (2\le i \le k).$ 

$\mathrm(ii)$\ \ $B=\bigcup_{i=1}^kI_i$, $I_i\cap I_j=\emptyset\ (i\ne j)$ and $\sharp B\ge n+1.$

$\mathrm(iii)$\ \ For each $1\le i\le k,$ there exist meromorphic functions $c_\alpha\in \mathcal {K}_f\setminus\{0\}$ such that 
$$\sum_{Q_\alpha (f)\in I_i}c_\alpha Q_\alpha (f)\in \left (\bigcup_{j=1}^{i-1} I_j\right )_{\mathcal K_f}.$$

\end{lemma}
\textbf{\textit{Proof.}}\  Denote by $V^d_f$ the vector space of all homogeneous polynomials of degree $d$ in $\mathcal K_f[x_0,...,x_n]$. It is seen that $\dim V^d_f=\binom{n+d}{n}=N+1$. 

$\bullet$ We set $A_0=\{Q_i(f)\ ;\ 0\le i\le n(N+1)\}$. We are going to construct the subset $B_0$ of $A_0$ as follows:

Since $\sharp A_0>N+1=\dim V^d_f$, the set $A_0$ is linearly independent over $\mathcal K_f$. Therefore, there exists a minimal subset $I^0_1$ over $\mathcal K_f$ of $A_0$. If $\sharp I^0_1\ge n+1$ or $\left (I^0_1\right )_{\mathcal K_f}\cap \left (A_0\setminus I^0_1\right )_{\mathcal K_f}=\{0\}$ then we stop the process and set $B_0=I^0_1, A_1=A_0\setminus B_0$. 

Otherwise, since $\left (I^0_1\right )_{\mathcal K_f}\cap \left (A_0\setminus I^0_1\right )_{\mathcal K_f}\ne\{0\}$, we now choose a subset $I^0_2$ of  $A_0\setminus I^0_1$ such that $I^0_2$ 
is the minimal subset of $A_0\setminus I^0_1$ satisfying $\left (I^0_1\right )_{\mathcal K_f}\cap \left (I^0_2\right )_{\mathcal K_f}\ne\{0\}$. By the minimality, the subset $I^0_2$ is linearly independent over $\mathcal K_f$. If $\sharp (I^0_1\cup I^0_2)\ge n+1$ or $\left (I^0_1\cup I^0_2\right )_{\mathcal K_f}\cap \left (A_0\setminus (I^0_1\cup I^0_2)\right )_{\mathcal K_f}=\{0\}$ then we stop the process and set $B_0=I^0_1\cup I^0_2, A_1=A_0\setminus B_0$. 

Otherwise, by repeating the above argument, we have a subset $I^0_3$ of $A_0\setminus (I^0_1\cup I^0_2).$

Continuiting this process, there exist subsets $I^0_1,...,I^0_k$ such that: $I^0_i$ is a subset of $A_0\setminus \bigcup_{j=1}^{i-1}I^0_j$, $I^0_j$ is linearly independent over $\mathcal K_f\ (2\le j\le k)$, $\left (I^0_i\right )_{\mathcal K_f}\cap \left (\bigcup_{j=1}^{i-1} I^0_j\right )_{\mathcal K_f}\ne\{0\}$, $\sharp B_0\ge n+1$ or $\left (B_0\right )_{\mathcal K_f}\cap \left (A_0\setminus B_0\right )_{\mathcal K_f}=\{0\}$. Also, by the minimality of each subset $I^0_i \ (2\le i\le k)$, there exist nonzero meromorphic functions $c^0_\alpha\in\mathcal K_f$ such that
$$ \sum_{Q_\alpha (f)\in I^0_i}c^0_\alpha Q_\alpha (f)\in\left (\bigcup_{j=1}^{i-1}I^0_j\right )_{\mathcal K_f}.$$

$\bullet$ If $\sharp B_0\ge n+1$, by setting $B=B_0, I_i=I^0_i$ then the proof is finished. 

Otherwise, we have $\left (B_0\right )_{\mathcal K_f}\cap \left (A_0\setminus B_0\right )_{\mathcal K_f}=\{0\}$. We set $A_1=A_0\setminus B_0$. Then $\dim (A_1)_{\mathcal K_f}\le N+1-\dim (B_0)_{\mathcal K_f}\le N$ and $\sharp A_1\ge nN+1> N\ge \dim (A_1)_{\mathcal K_f}$. Similarly, we construct the subset $B_1$ of $A_1$ with the same properties as $B_0$. 

$\bullet$ If $\sharp B_1\ge n+1$ then the proof is finished. Otherwise, by repeating the same argument we have subsets $A_3,B_3$ and $I^3_i$.

Continuiting this process, we have the following two cases:

\textbf{Case 1.} By this way, we may construct subsets $B_1,...,B_{N}$ with $\sharp B_i\le n\ (1\le i\le N)$. We set $B_{N+1}=A_0\setminus\bigcup_{i=0}^{N}B_i$. Then $\sharp B_{N+1}\ge n(N+1)+1-n(N+1)=1$. Then $\dim \left (B_{N+1}\right )_{\mathcal K_f}\ge 1$. On the other hand, it is easy to see that
$$ \dim \left (B_{N+1}\right )_{\mathcal K_f}=\dim \left (A_0\right )_{\mathcal K_f}-\sum_{i=0}^N\dim \left (B_i\right )_{\mathcal K_f}\le N+1-(N+1)=0.$$
This is a contradiction. Hence this case is impossible.

\textbf{Case 2.} At the step $k-th$ $(k\le N)$, we get $\sharp B_k\ge n+1$. Then similarly as above, the proof is finished.\hfill$\square$

\begin{lemma}\label{3.2} 
Let $f$ be as in Theorem \ref{1.1}. Let $\{Q_i\}_{i=0}^{n(N+1)}$ be a set of  homogeneous polynomials in $\mathcal K_f[x_0,...,x_n]$ of common degree $d$ in weakly general position, where $N=\binom{n+d}{n}-1$. Assume that $Q_i(f)\not\equiv 0\ (0\le i\le n(N+1))$. Then we have
$$ ||\ T_f(r)\le \sum_{i=0}^{n(N+1)}\dfrac{1}{d}N^{[N]}_{Q_i(f)}(r)+o(T_f(r)).$$
\end{lemma}
\textbf{\textit{Proof.}}\
By Lemma \ref{3.1}, we may assume that there exist subsets 
$$I_i=\{Q_{t_{i}+1}(f) ,..., Q_{t_{i+1}}(f)\}\ \ (1\le i \le k)$$
 and functions $c_i\in\mathcal K_f\setminus\{0\}\ (t_2+1\le i\le t_{k+1}),$
where $t_1=-1,$ which satisfy the assertions of Lemma \ref{3.1}.

Since $I_1$ is minimal over $\mathcal K_f$, there exist $c_{1j}\in\mathcal {R}\setminus \{0\}$ such that 
$$\sum_{j=0}^{t_2}c_{1j} Q_j(f)=0.$$

Define $c_{1j}=0$ for all $j>t_1.$ Then $\sum_{j=0}^{t_{k+1}} c_{1j} Q_j(f)=0.$

Since $\{c_{1j}Q_j(f)\}_{j=1}^{t_2}$ is linearly independent over $\mathcal K_f,$ there exists
an admissible set $\{\alpha_{11},...,\alpha_{1t_2}\}\subset \Z^m_+$ \ $(|\alpha_{1j}|\le t_2-1\le N)$ such that 
\begin{align*} 
\ A_1\equiv &\left | \begin {array}{cccc}
\mathcal {D}^{\alpha_{11}}(c_{11}Q_1(f)) &\cdots & \mathcal {D}^{\alpha_{11}}(c_{1t_2}Q_{t_2}(f)) \\
\mathcal {D}^{\alpha_{12}}(c_{11}Q_1(f)) &\cdots & \mathcal {D}^{\alpha_{12}}(c_{1t_2}Q_{t_2}(f)) \\
\vdots &\vdots &\vdots\\
\mathcal {D}^{\alpha_{1t_2}}(c_{11}Q_1(f))&\cdots &\mathcal {D}^{\alpha_{1t_2}}(c_{1t_2}Q_{t_2}(f))\\
\end {array}
\right|\\\\
\equiv &f_0^{t_1}\cdot \left | \begin {array}{cccc}
\mathcal {D}^{\alpha_{11}}\biggl(\dfrac {c_{11}Q_1(f)}{Q_0(f)}\biggl) &\cdots & \mathcal {D}^{\alpha_{11}}\biggl(\dfrac{c_{1t_2}Q_{t_2}(f)}{Q_0(f)}\biggl) \\
\mathcal {D}^{\alpha_{12}}\biggl(\dfrac {c_{11}Q_1(f))}{Q_0(f)}\biggl) &\cdots & \mathcal {D}^{\alpha_{12}}\biggl(\dfrac{c_{1t_2}Q_{t_2}(f)}{Q_0(f)}\biggl) \\
\vdots &\vdots &\vdots\\
\mathcal {D}^{\alpha_{1t_2}}\biggl(\dfrac{c_{11}Q_1(f)}{Q_0(f)}\biggl)&\cdots &\mathcal {D}^{\alpha_{1t_1}}\biggl(\dfrac{c_{1t_2}Q_{t_2}(f)}{Q_0(f)}\biggl)\\
\end {array}
\right|\equiv (Q_0(f))^{t_2}\cdot \tilde A_1\not\equiv 0.
\end{align*}

Now consider $i\ge 2.$ We set $c_{ij}=c_j\not\equiv 0\ (t_{i}+1\le j\le t_{i+1})$, then $\sum_{j=t_{i}+1}^{t_{i+1}}c_{ij}Q_j(f)\in \left(\bigcup_{j=1}^{i-1}I_j \right)_{\mathcal K_f}.$ 
Therefore, there exist meromorphic functions $c_{ij}\in \mathcal {K}_f\ (0\le j\le t_i)$ 
such that $\sum_{j=0}^{t_{i+1}}c_{ij}Q_j(f)=0.$

Define $c_{ij}=0$ for all $j>t_{i+1}.$ Then $\sum_{j=0}^{t_{k+1}}c_{ij}Q_j(f)=0.$

Since $\{c_{ij}Q_j(f)\}_{j=t_{i}+1}^{t_{i+1}}$ is linearly independent over $\mathcal K_f$, 
there exists $\{\alpha_{ij}\}_{j=t_{i}+1}^{t_{i+1}}\subset \Z^m_+$ \ $(|\alpha_{ij}|\le t_{i+1}-t_{i}-1\le N)$ such that
\begin{align*}
A_i=&\det \biggl(\mathcal{D}^{\alpha_{ij}}\biggl(c_{is}Q_s(f)\biggl)\biggl)_{j,s=t_{i}+1}^{t_{i+1}}=(Q_0(f))^{t_{i+1}-t_i}\cdot 
\det\biggl(\mathcal D^{\alpha_{ij}}\biggl(\dfrac{c_{is}Q_s(f)}{Q_0(f)}\biggl)\biggl)_{j,s=t_{i}+1}^{t_{i+1}}\\
=&Q_0(f)^{t_{i+1}-t_i}\cdot \tilde A_i\not\equiv 0.
\end{align*}
Consider an $t_{k+1}\times (t_{k+1}+1)$ minor matrixes  $\mathcal {T}$ and $\mathcal{\tilde T}$ given by
\begin{align*}
\mathcal {T}=&\left [ \begin {array}{cccc}
\mathcal {D}^{\alpha_{11}}(c_{10}Q_0(f)) &\cdots & \mathcal {D}^{\alpha_{11}}(c_{1t_{k+1}}Q_{t_{k+1}}(f)) \\
\mathcal {D}^{\alpha_{12}}(c_{10}Q_0(f)) &\cdots & \mathcal {D}^{\alpha_{12}}(c_{1t_{k+1}}Q_{t_{k+1}}(f)) \\
\vdots &\vdots &\vdots\\
\mathcal {D}^{\alpha_{1t_2}}(c_{10}Q_0(f))&\cdots &\mathcal {D}^{\alpha_{1t_2}}(c_{1t_{k+1}}Q_{t_{k+1}}(f))\\
\mathcal {D}^{\alpha_{2t_2+1}}(c_{20}Q_0(f)) &\cdots & \mathcal {D}^{\alpha_{2t_2+1}}(c_{2t_{k+1}}Q_{t_{k+1}}(f)) \\
\mathcal {D}^{\alpha_{2t_2+2}}(c_{20}Q_0(f)) &\cdots & \mathcal {D}^{\alpha_{2t_2+2}}(c_{2t_{k+1}}Q_{t_{k+1}}(f)) \\
\vdots &\vdots &\vdots\\
\mathcal {D}^{\alpha_{2t_3}}(c_{20}Q_0(f))&\cdots &\mathcal {D}^{\alpha_{2t_3}}(c_{2t_{k+1}}Q_{t_{k+1}}(f))\\
\vdots &\vdots &\vdots\\
\mathcal {D}^{\alpha_{kt_{k}+1}}(c_{k0}Q_0(f)) &\cdots & \mathcal {D}^{\alpha_{kt_{k}+1}}(c_{kt_{k+1}}Q_{t_{k+1}}(f)) \\
\mathcal {D}^{\alpha_{kt_{k}+2}}(c_{k0}Q_0(f)) &\cdots & \mathcal {D}^{\alpha_{kt_{k}+2}}(c_{kt_{k+1}}Q_{t_{k+1}}(f)) \\
\vdots &\vdots &\vdots\\
\mathcal {D}^{\alpha_{kt_{k+1}}}(c_{k0}Q_0(f))&\cdots &\mathcal {D}^{\alpha_{kt_{k+1}}}(c_{kt_{k+1}}Q_{t_{k+1}}(f))
\end {array}
\right]\\\\
\mathcal {\tilde T}=&
\left [ \begin {array}{cccc}
\mathcal {D}^{\alpha_{11}}\biggl(\dfrac{c_{10}Q_0(f)}{Q_0(f)}\biggl) &\cdots & \mathcal {D}^{\alpha_{11}}\biggl(\dfrac{c_{1t_{k+1}}Q_{t_{k+1}}(f)}{Q_0(f)}\biggl) \\
\vdots &\vdots &\vdots\\
\mathcal {D}^{\alpha_{1t_2}}\biggl(\dfrac{c_{10}Q_0(f)}{Q_0(f)}\biggl)&\cdots &\mathcal {D}^{\alpha_{1t_2}}\biggl(\dfrac{c_{1t_{k+1}}Q_{t_{k+1}}(f)}{Q_0(f)}\biggl)\\
\mathcal {D}^{\alpha_{2t_2+1}}\biggl(\dfrac{c_{20}Q_0(f)}{Q_0(f)}\biggl) &\cdots & \mathcal {D}^{\alpha_{2t_2+1}}\biggl(\dfrac{c_{1t_{k+1}}Q_{t_{k+1}}(f)}{Q_0(f)}\biggl) \\
\vdots &\vdots &\vdots\\
\mathcal {D}^{\alpha_{2t_3}}\biggl(\dfrac{c_{20}Q_0(f)}{Q_0(f)}\biggl)&\cdots &\mathcal {D}^{\alpha_{2t_3}}\biggl(\dfrac{c_{2t_{k+1}}Q_{t_{k+1}}(f)}{Q_0(f)}\biggl)\\
\vdots &\vdots &\vdots\\
\mathcal {D}^{\alpha_{kt_{k}+1}}\biggl(\dfrac{c_{k0}Q_0(f)}{Q_0(f)}\biggl) &\cdots & \mathcal {D}^{\alpha_{kt_{k}+1}}\biggl(\dfrac{c_{kt_{k+1}}Q_{t_{k+1}}(f)}{Q_0(f)}\biggl) \\
\vdots &\vdots &\vdots\\
\mathcal {D}^{\alpha_{kt_{k+1}}}\biggl(\dfrac{c_{k0}Q_0(f)}{Q_0(f)}\biggl)&\cdots &\mathcal {D}^{\alpha_{kt_{k+1}}}\biggl(\dfrac{c_{kt_{k+1}}Q_{t_{k+1}}(f)}{Q_0(f)}\biggl)
\end {array}
\right].
\end{align*}
Denote by $D_i$ (resp. ${\tilde D}_i$) the determinant of the matrix obtained  by deleting the $(i+1)$-th column 
of the minor  matrix  $\mathcal T$ (resp. $\mathcal{\tilde T}$). It is clear that the sum of each row of \ $\mathcal T$ (resp.$\mathcal{\tilde T}$) is zero, then we have
\begin{align*}
D_i&={(-1)}^iD_0={(-1)}^i\prod_{i=1}^{k}A_i={(-1)}^i(Q_0(f))^{t_{k+1}}\prod_{i=1}^{k}\tilde A_i\\
&={(-1)}^i(Q_0(f))^{t_{k+1}}{\tilde D}_0=(Q_0(f))^{t_{k+1}}{\tilde D}_i .
\end{align*}

Since $\sharp (\bigcup_{i=1}^{k}I_i)\ge n+1$ and $Q_0,...,Q_{t_{k+1}}$ are in weakly general position, by Lemma 2.8 there exists a function $\Psi\in \mathcal C_f$ such that
$$\ \ ||f(z)||^d\le \Psi (z)\cdot \max_{0\le i\le t_{k+1}}\bigl (|Q_i(f)(z)|\bigl )\ (z\in \C^m).$$

Fix $z_0 \in \C^m.$ Take $i\ (0\le i \le t_k)$ such that $|Q_i(f)(z_0)|=\max_{0\le j\le t_k}|Q_j(f)(z_0)|.$ 
Then
$$
\dfrac{|D_0(z_0)|\cdot ||f(z_0)||^d}{\prod_{j=0}^{t_{k+1}}|Q_j(f)(z_0)|}
=\dfrac{|D_i(z_0)|}{\prod_{\underset{j\ne i}{j=0}}^{t_{k+1}}|Q_j(f)(z_0)|}\cdot \biggl(\dfrac
{||f(z_0)||^d}{|Q_i(f)(z_0)|}\biggl)\le \Psi (z_0)\cdot \dfrac{|D_i(z_0)|}{\prod_{\underset{j\ne i}{j=0}}^{t_{k+1}}|Q_j(f)(z_0)|}.
$$
This implies that
\begin{align*}
\log\dfrac{|D_0(z_0)|.||f(z_0)||^d}{\prod_{j=0}^{t_{k+1}}|Q_j(f)(z_0)|}& \le\log^+\biggl ( \Psi (z_0)\cdot \biggl(\dfrac{|D_i(z_0)|}{\prod_{j=0,j\ne i}^{t_{k+1}}|Q_j(f)(z_0)|}\biggl)\biggl  )\\
& \le\log^+\biggl(\dfrac{|D_i(z_0)|}{\prod_{j=0,j\ne i}^{t_k}|Q_j(f)(z_0)|}\biggl)+\log^+\Psi (z_0).
\end{align*}
Thus, for each $z\in \C^m,$ we have
\begin{align}\nonumber
\log\dfrac{|D_0(z)|.||f(z)||^d}{\prod_{i=0}^{t_{k+1}}|Q_i(f)(z)|}& \le\sum_{i=0}^{t_{k+1}}\log^+\biggl(\dfrac{|D_i(z)|}{\prod_{j=0,j\ne i}^{t_k}|Q_j(f)(z)|}\biggl)+\log^+ \Psi (z)\\
\label{3.3}
& =\sum_{i=0}^{t_{k+1}}\log^+\biggl(\dfrac{|{\tilde D}_i(z)|}{\prod_{j=0,j\ne i}^{t_k}\biggl|\dfrac{Q_j(f)(z)}{Q_0(f)(z)}\biggl|}\biggl)+\log^+ \Psi (z).
\end{align}
Note that \  
$\dfrac{{\tilde D}_i}{\prod_{j=0,j\ne i}^{t_{k+1}}\dfrac{Q_j(f)}{Q_0(f)}}=\det$
$\left [ \begin {array}{cccc}
\dfrac {\mathcal {D}^{\alpha_{11}}\biggl(\dfrac {c_{10}Q_0(f)}{Q_0(f)}\biggl)}{\dfrac{Q_0(f)}{Q_0(f)}} &\cdots &
\dfrac {\mathcal {D}^{\alpha_{11}}\biggl(\dfrac {c_{1t_{k+1}}Q_{t_{k+1}}(f)}{Q_0(f)}\biggl)}{\dfrac{Q_{t_{k+1}}(f)}{Q_0(f)}} \\
\vdots &\vdots &\vdots \\
\dfrac {\mathcal {D}^{\alpha_{kt_{k+1}}}\biggl(\dfrac {c_{k0}Q_0(f)}{Q_0(f)}\biggl)}{\dfrac{Q_0(f)}{Q_0(f)}} &\cdots &
\dfrac {\mathcal {D}^{\alpha_{kt_{k+1}}}\biggl(\dfrac {c_{kt_{k+1}}Q_{t_{k+1}}(f)}{Q_0(f)}\biggl)}{\dfrac {Q_{t_{k+1}}(f)}{Q_0(f)}} \\
\end {array}
\right]$

\quad (The determinant is counted after deleting the $i$-th column in the above matrix)

By the lemma on logarithmic derivative, for each $i$ and $c\in\mathcal K_f$ we have
\begin{align*}
\biggl|\biggl| \quad\quad m \biggl(r,\dfrac {\mathcal {D}^{\alpha}\biggl(\dfrac {cQ_j(f)}{Q_0(f)}\biggl)}{\dfrac{Q_j(f)}{Q_0(f)}}\biggl)&\le
m \biggl(r,\dfrac {\mathcal {D}^{\alpha}\biggl(\dfrac {cQ_j(f)}{Q_0(f)}\biggl)}{\dfrac{cQ_j(f)}{Q_0(f)}}\biggl)+m(r,c)\\
&\le O\biggl(\log^+T_{\dfrac {cQ_j(f)}{Q_0(f)}}(r)\biggl )+T_c(r)=o(T_f(r))
\end{align*}
Therefore, we have
$$\biggl|\biggl| \quad m\biggl(r,\dfrac{{\tilde D}_i}{\prod_{j=0,j\ne i}^{t_{k+1}}\dfrac{Q_j(f)}{Q_0(f)}}\biggl)=o(T_f(r))\ (0 \le i \le t_k).$$

Integrating both sides of the inequality (\ref{3.3}), we get
\begin{align*}
\biggl|\biggl|  \  \int_{S(r)}\log ||f||^d \sigma_m &+ \int_{S(r)}\log \biggl(\dfrac{|D_0|}{\prod_{i=0}^{t_{k+1}} |Q_i(f)|} \biggl)\sigma_m\\
&\le \sum_{i=0}^{t_{k+1}} \int_{S(r)}\log^+ \biggl(\dfrac{|{\tilde D}_i|}{\prod_{j=0,j\ne i}^{t_{k+1}}
|\dfrac{Q_j(f)}{Q_0(f)}|}\biggl)\sigma_m +\int_{S(r)}\log^+ \Psi(z)\sigma_m\\
&\le \sum_{i=0}^{t_{k+1}} m\biggl(r,\dfrac{{\tilde D}_i}{\prod_{j=0,j\ne i}^{t_{k+1}}\dfrac{Q_j(f)}{Q_0(f)}}\biggl)+o(T_f(r))=o(T_f(r)).
\end{align*}
By Jensen formula, the above inequality implies that
\begin{align}\label{3.4}
||\ \ dT_f(r)+N_{D_0}(r)-N_{\frac{1}{D_0}}(r)-\sum_{i=0}^{t_{k+1}}N_{Q_i(f)}(r)\le o(T_f(r)).
\end{align}
We see that a pole of $D_0$ must be pole of some $c_{is}$ or pole of some nonzero coefficients $a_{iI}$ of  $Q_i$ and
$$ N_{\frac{1}{D_0}}(r)\le O(\sum_{i,s}N_{\frac{1}{c_{is}}}(r)+\sum_{a_{iI}\not\equiv 0}N_{\frac{1}{a_{iI}}}(r)) =o(T_f(r)).$$
Therefore, the inequality (\ref{3.4}) implies that
\begin{align}\label{3.5}
||\ \ dT_f(r)\le \sum_{i=0}^{t_{k+1}}N_{Q_i(f)}(r)-N_{D_0}(r)+o(T_f(r)).
\end{align}
Here we note that $D_i=(-1)^iD_0$, then $\nu^0_{D_i}=\nu^0_{D_0}$.

We now assume that $z$ is a zero of some functions $Q_i(f)$. Since $t_{k+1}+1\ge n+1$ and $z$ can not be zero of more than $n$ functions $Q_i(f)$,  without loss of generality we may assume that $z$ is not zero of $Q_0(f)$. Then
\begin{align*}
&\nu^0_{\mathcal {D}^{\alpha_{st_{s-1}+j}}(c_{si}Q_i(f))}(z)
\ge \min_{\beta\in \Z_+^m \text { with } \alpha_{st_{s-1}+j}-\beta \in\Z_+^m}
\{\nu^0_{\mathcal{D}^{\beta}c_{si}\mathcal D^{\alpha_{st_{s-1}+j}-\beta }Q_i(f)}(z)\}\\
&\ge\min_{\beta\in \Z_+^m \text { with } \alpha_{st_{s-1}+j}-\beta \in\Z_+^m}\bigl{\{}\max\{0,\nu^0_{Q_i(f)}(z)-|\alpha_{st_{s-1}+j}-\beta|\}-(\beta+1)\nu^{\infty}_{c_{si}}(z)\bigl{\}}\\
&\ge\max\{0,\nu_{Q_i(f)}^{0}(z)-N\}-(N+1)\nu_{c_{si}}^{\infty}(z)
\end{align*}
for each $1 \le i \le t_{k+1}, 1\le j \le t_s-t_{s-1}, 1\le s \le k+1,$ where $t_0=0.$. 

Put
$I(z)=(N+1)\sum_{s=1}^{k}\sum_{i=0}^{t_k}(t_{s}-t_{s-1})\nu_{c_{si}}^{\infty}(z).$
Then 
\begin{align}\label{3.6}
\nu_{D_0}(z)\ge\sum_{i=0}^{t_{k+1}}\max\{0,\nu_{Q_i(f)}^{0}(z)-N\}-I(z).
\end{align}
We note that if $z$ is not zero of a function $Q_i(f)$ with $i\ne 0$, replacing $D_0$ by $D_i$ and repeating the same above argument we again get the inequality (\ref{3.6}). Hence (\ref{3.6}) holds for all $z\in\C^m$. It follows that
$$ \sum_{i=0}^{t_{k+1}}\nu_{Q_i(f)}^{0}(z)-\nu_{D_0}(z)\le \sum_{i=0}^{t_k-1}\min\{N,\nu_{Q_i(f)}^{0}(z)\}+I(z). $$
Integrating both sides of the above inequality, we get
$$ \sum_{i=0}^{t_{k+1}}N_{Q_i(f)}(r)-N_{D_0}(r)\le \sum_{i=0}^{t_{k+1}}N^{[N]}_{Q_i(f)}(r)+o(T_f(r)).$$
Combining this and (\ref{3.5}), we get
$$ ||\ T_f(r)\le \sum_{i=0}^{n(N+1)}\dfrac{1}{d}N^{[N]}_{Q_i(f)}(r)+o(T_f(r)).$$
The lemma is proved.\hfill$\square$
\vskip0.2cm
\noindent
{\bf Proof of Theorem \ref{1.1}.}

\vskip0.2cm
We first prove the theorem for the case where all $Q_i\ (i=1,...,q)$ have the same degree $d$. By changing the homogeneous coordinates of $\P^n(\C)$ if necessary, we may assume that $a_{iI_1}\not\equiv 0$ for every $i=1,...,q$. We set $\tilde Q_i=\dfrac{1}{a_{iI_1}}Q_i$. Then $\{\tilde Q_i\}_{i=1}^q$ is a set of homogeneous polynomials in $\mathcal K_f [x_0,...,x_n]$ in weakly general position.

Consider $(nN+n+1)$  polynomials $\tilde Q_{i_1},...,\tilde Q_{i_{nN+n+1}}\ (1\le i_j\le q).$ Applying Lemma \ref{3.2}, we have
$$\bigl|\bigl|\quad T_f(r)\le \sum_{j=1}^{nN+n+1}\dfrac{1}{d}N^{[N]}_{\tilde Q_i(f)}(r)+o(T_f(r))\
\le \sum_{j=1}^{nN+n+1}\dfrac{1}{d}N^{[N]}_{Q_i(f)}(r)+o(T_f(r)).$$
Taking summing-up of both sides of this inequality over all combinations $\{i_1,...,i_{nN+n+1}\}$ with 
$1\le i_1< ...<i_{nN+n+1}\le q,$ we have
$$\biggl|\biggl|\quad\dfrac{q}{nN+n+1}T_f(r)\le \sum_{j=1}^{nN+n+1}\dfrac{1}{d}N^{[N]}_{Q_i(f)}(r)+o(T_f(r)).$$
The theorem is proved in this case.

We now prove the theorem for the general case where $\deg Q_i=d_i$. Then, applying the above case for $f$ and the moving hypersurfaces $Q^{\frac{d}{d_i}}_i\ (i=1,...,q)$ of common degree $d$, we have
\begin{align*}
\biggl|\biggl|\quad\dfrac{q}{nN+n+1}T_f(r)&\le \sum_{j=1}^{q}\dfrac{1}{d}N^{[N]}_{Q^{d/d_i}_i(f)}(r)+o(T_f(r))\\
&\le \sum_{j=1}^{q}\dfrac{1}{d}\frac{d}{d_i}N^{[N]}_{Q_i(f)}(r)+o(T_f(r))\\
&=\sum_{j=1}^{q}\dfrac{1}{d_i}N^{[N]}_{Q_i(f)}(r)+o(T_f(r)).
\end{align*}
The theorem is proved.\hfill$\square$

\noindent
{\bf Proof of Theorem \ref{1.2}.}

By repeating the argument as in the proof of Theorem \ref{1.1}, it suffices to prove the theorem for the case where all  $Q_i$ have the same degree.

By changing the homogeneous coordinates of $\P^n(\C)$ if necessary, we may assume that $a_{iI_1}\not\equiv 0$ for every $i=1,...,q$. We set $\tilde Q_i=\dfrac{1}{a_{iI_1}}Q_i$. Then $\{\tilde Q_i\}_{i=1}^q$ is a set of homogeneous polynomials in $\mathcal K_f [x_0,...,x_n]$ in weakly general position.

Consider $(N+2)$  polynomials $\tilde Q_{i_1},...,\tilde Q_{i_{N+2}}\ (1\le i_j\le q)$. We see that $\dim (\tilde Q_{i_j}\ ;\ 1\le j\le N+2)_{\mathcal{\tilde K}_{\{Q_i\}_{i=1}^q}}\le N+1<N+2$. Then the set $\{Q_{i_1},...,Q_{i_{N+2}}\}$ is linearly independent over $\mathcal{\tilde K}_{\{Q_i\}_{i=1}^q}$. Hence, there exists a minimal subset over $\mathcal{\tilde K}_{\{Q_i\}_{i=1}^q}$, for instance that is $\{\tilde Q_{i_1},...,\tilde Q_{i_{t}}\}$, of  $\{\tilde Q_{i_1},...,\tilde Q_{i_{N+2}}\}$. Then, there exist nonzero functions 
$c_j\ (1\le j\le t)$ in $\mathcal{\tilde K}_{\{Q_i\}_{i=1}^q}$ such that
$$ c_1\tilde Q_{i_1}+\cdots +c_t\tilde Q_{i_t}=0. $$
Since $Q_{i_1},....,Q_{i_{N+2}}$ are in weakly general position, $t\ge n+2.$ Setting $F_j=c_jQ_j(f)$, we have
$$ F_1+\cdots F_{t-1}=-F_t. $$
Choose a meromorphic functions $h$ so that $F=(hF_1:\cdots :hF_{t-1})$ is a reduced representation of a meromorphic mapping $F$ from $\C^m$ into $\P^n(\C)$. It is seen that 
$$ N_h(r)\le \sum_{j=1}^{t-1}(N_{\frac{1}{c_j}}(r)+N_{a_{i_jI_1}}(r) )=o(T_f(r)).$$
On the other hand, by the minimality of the set $\{\tilde Q_{i_1},...,\tilde Q_{i_{t}}\}$, then $F$ is linearly nondegenerate over $\C$. Applying the second main theorem for fixed hyperplanes, we get
\begin{align*}
||\ T_F(r)&\le \sum_{j=1}^tN_{hF_j}^{[t-2]}(r)+o(T_F(r))\\
&\le \sum_{j=1}^t(N_{\tilde Q_{i_j}(f)}^{[t-2]}(r)+N_{c_j}^{[t-2]}(r))+tN_{h}^{[t-2]}(r)+o(T_F(r))\\
&=\sum_{j=1}^tN_{Q_{i_j}(f)}^{[t-2]}(r)+o(T_f(r))\le \sum_{j=1}^{N+2}N_{Q_{i_j}(f)}^{[N]}(r)+o(T_f(r)).
\end{align*}
It follows that
\begin{align*}
||\ T_f(r)=\dfrac{1}{d}T_F(r)+o(T_f(r))\le \sum_{j=1}^{N+2}\dfrac{1}{d}N_{Q_{i_j}(f)}^{[N]}(r)+o(T_f(r)).
\end{align*}
Taking summing-up of both sides of this inequality over all combinations $\{i_1,...,i_{N+2}\}$ with 
$1\le i_1< ...<i_{N+2}\le q,$ we have
$$\biggl|\biggl|\quad\dfrac{q}{N+2}T_f(r)\le \sum_{j=1}^{q}\dfrac{1}{d}N^{[N]}_{Q_i(f)}(r)+o(T_f(r)).$$
The theorem is proved.\hfill$\square$

\section{Uniqueness problem of meromorphic mappings sharing moving hypersurfaces}

In order to prove Theorem \ref{1.3} and Theorem \ref{1.4} we need the following.
\begin{lemma}\label{4.1}
Let $f$ and $g$ be nonconstant meromorphic mappings of $\C^m$ into $\P^n(\C )$. Let $Q_i\ (i=1,...,q)$ be slow (with respect to $f$ and $g$) moving hypersurfaces in $\P^n(\C)$ in weakly general position with $\deg Q_i=d_i$. Assume that $\min\{\nu^0_{Q_i(f)}(z),1\}=\min\{\nu^0_{Q_i(g)}(z),1\}$ for all $1\le i\le q$. Put $d=lcm (d_1,...,d_q)$ and $N=\binom{n+d}{n}-1$. Then the following assertions hold:

(i) If $q>\frac{2N(nN+n+1)}{d}$ then $ ||\ T_f(r)=O(T_g(r))\ \text{ and }\  ||\ T_g(r)=O(T_f(r)).$

(ii) If both $f$ and $g$ are algebraically nondegenerate over $\mathcal{\tilde K}_{\{Q_i\}_{i=1}^q}$ and $q\ge n+2$ then $ ||\ T_f(r)=O(T_g(r))\ \text{ and }\  ||\ T_g(r)=O(T_f(r)).$
\end{lemma}
\textbf{\textit{Proof.}}\
(i)\  It is clear that $q>nN+n+1.$ Then applying Theorem \ref{1.1} for $f$, we have

\begin{align*}
||\  \dfrac{q}{nN+n+1}T_g(r)\le &\sum_{i=1}^q \dfrac{1}{d_i}N_{Q_i(g)}^{[N]}(r)+o(T_g(r))\\
\le &\sum_{i=1}^q\dfrac{N}{d_i}\ N_{Q_i(g)}^{[1]}(r)+o(T_g(r))\\
\le &\sum_{i=1}^q\dfrac{N}{d_i}\ N_{Q_i(f)}^{[1]}(r)+o(T_g(r))\\
\le & qN\ T_f(r)+o(T_g(r)).
\end{align*}
Hence \quad $|| \quad T_g(r)=O(T_f(r)).$ Similarly, we get \  \ $|| \ \ T_f(r)=O(T_g(r)).$

(ii) Applying Theorem A with $\epsilon=\dfrac{1}{2}$, then there exists a positive integer $L$ such that
\begin{align*}
||\ (q-n-\dfrac{3}{2})T_f(r)\le &\sum_{i=1}^q\dfrac{1}{d_i}N^{[L]}_{Q_i(f)}(r)+o(T_f(r)),\\ 
||\ (q-n-\dfrac{3}{2})T_g(r)\le &\sum_{i=1}^q\dfrac{1}{d_i}N^{[L]}_{Q_i(g)}(r)+o(T_g(r)).
\end{align*}
Therefore, we have
\begin{align*}
||\  (q-n-\dfrac{3}{2})T_f(r)\le &\sum_{i=1}^q\dfrac{1}{d_i}N^{[L]}_{Q_i(f)}(r)+o(T_f(r))\le \sum_{i=1}^q\dfrac{L}{d_i}\ N^{[1]}_{Q_i(f)}(r)+o(T_f(r))\\ 
=  &\sum_{i=1}^q\dfrac{L}{d_i}\ N_{Q_i(g)}^{[1]}(r)+o(T_f(r))\le qL\ T_g(r)+o(T_g(r)).
\end{align*}
Hence \quad $|| \quad T_f(r)=O(T_g(r)).$ Similarly, we get \  \ $|| \ \ T_g(r)=O(T_f(r)).$

\vskip0.2cm
\noindent
\textbf{Proof of Theorem \ref{1.3}.}\ We assume that $f$ and $g$ have reduced representations $f=(f_0:\cdots :f_n)$ and $g=(g_0:\cdots :g_n)$ respectively.

a)\ By Lemma \ref{4.1} (i) , we have $ ||\ T_f(r)=O(T_g(r))\ \text{ and }\  ||\ T_g(r)=O(T_f(r)).$ Suppose that $f$ and $g$ are two distinct maps. Then there exist two index $s,t \ (0\le s<t\le n)$ satisfying
$$ H:=f_sg_t-f_tg_s\not\equiv 0. $$
Set $S=\bigcup\{\bigcap_{j=0}^k\zero Q_{i_j}(f)\ ;\ 1\le i_0<\cdots <i_k\le q\}$. Then $S$ is either an analytic subset of codimension at least two of $\C^m$ or an empty set. 

Assume that $z$ is a zero of some $Q_i(f)\ (1\le i\le q)$ and $z\not\in S$. Then the condition (iii) yields that $z$ is a zero of the function $H.$ Also, since $z\not\in S$, $z$ can not be zero of more than $k$ functions $Q_i(f)$. Therefore, we have
$$ \nu^0_H(z) =1\ge \dfrac{1}{k}\sum_{i=1}^q\min\{1,\nu^0_{Q_i(f)}(z)\}.$$
This inequality holds for every $z$ outside the analytic subset $S$ of codimension at least two. Then, it follows that
\begin{align}\label{4.2}
N_H(r)\ge \dfrac{1}{k}\sum_{i=1}^qN^{[1]}_{Q_i(f)}(r).
\end{align}

On the other hand, by the definition of the characteristic function and Jensen formula, we have
\begin{align*}
N_H(r)&=\int_{S(r)}\log |f_sg_t-f_tg_s|\sigma_m\\ 
& \le \int_{S(r)}\log ||f||\sigma_m +\int_{S(r)}\log ||f||\sigma_m \\
&=T_f(r)+T_g(r).
\end{align*}
Combining this and (\ref{4.2}), we obtain
$$ T_f(r)+T_g(r) \ge \dfrac{1}{k}\sum_{i=1}^qN^{[1]}_{Q_i(f)}(r).$$
Similarly, we have
$$ T_f(r)+T_g(r) \ge \dfrac{1}{k}\sum_{i=1}^qN^{[1]}_{Q_i(g)}(r).$$
Summing-up both sides of the above two inequalities, we have
\begin{align}
\nonumber
 2(T_f(r)+T_g(r))& \ge \dfrac{1}{k}\sum_{i=1}^qN^{[1]}_{Q_i(f)}(r)+\dfrac{1}{k}\sum_{i=1}^qN^{[1]}_{Q_i(g)}(r)\\
\nonumber
&=\dfrac{1}{k}\sum_{i=1}^qN^{[1]}_{Q^{d/d_i}_i(f)}(r)+\dfrac{1}{k}\sum_{i=1}^qN^{[1]}_{Q^{d/d_i}_i(g)}(r)\\
\label{4.3}
&\ge \sum_{i=1}^q\dfrac{1}{kN}N^{[N]}_{Q^{d/d_i}_i(f)}(r)+\sum_{i=1}^q\dfrac{1}{kN}N^{[N]}_{Q^{d/d_i}_i(g)}(r).
\end{align}
From (\ref{4.3}) and applying Theorem \ref{1.1} for $f$ and $g$, we have
\begin{align*}
2(T_f(r)+T_g(r))&\ge \sum_{i=1}^q\dfrac{1}{kN}N^{[N]}_{Q^{d/d_i}_i(f)}(r)+\sum_{i=1}^q\dfrac{1}{kN}N^{[N]}_{Q^{d/d_i}_i(g)}(r)\\
&\ge \dfrac{d}{kN}\dfrac{q}{nN+n+1}(T_f(r)+T_g(r))+o(T_f(r)+T_g(r)).
\end{align*}
Letting $r\longrightarrow +\infty$, we get $2\ge \frac{d}{kN}\frac{q}{nN+n+1}\Leftrightarrow q\le \frac{2kN(nN+n+1)}{d}.$ This is a contradiction.

Hence $f=g$. The assertion a) is proved.

b) By Lemma \ref{4.1} (ii) , we have $ ||\ T_f(r)=O(T_g(r))\ \text{ and }\  ||\ T_g(r)=O(T_f(r)).$ Suppose that $f$ and $g$ are two distinct maps. Repeating the same argument as in a), we get the following inequality, which is similar to (\ref{4.3}),
\begin{align}\label{4.4}
 2(T_f(r)+T_g(r))\ge \sum_{i=1}^q\dfrac{1}{kN}N^{[N]}_{Q^{d/d_i}_i(f)}(r)+\sum_{i=1}^q\dfrac{1}{kN}N^{[N]}_{Q^{d/d_i}_i(g)}(r).
\end{align}
From (\ref{4.4}) and applying Theorem \ref{1.2} for $f$ and $g$, we have
\begin{align*}
2(T_f(r)+T_g(r))&\ge \sum_{i=1}^q\dfrac{1}{kN}N^{[N]}_{Q^{d/d_i}_i(f)}(r)+\sum_{i=1}^q\dfrac{1}{kN}N^{[N]}_{Q^{d/d_i}_i(g)}(r)\\
&\ge \dfrac{d}{kN}\dfrac{q}{N+2}(T_f(r)+T_g(r))+o(T_f(r)+T_g(r)).
\end{align*}
Letting $r\longrightarrow +\infty$, we get $2\ge \frac{d}{kN}\frac{q}{N+2}\Leftrightarrow q\le \frac{2kN(N+2)}{d}.$ This is a contradiction.

Hence $f=g$. The assertion b) is proved.\hfill$\square$

\vskip0.2cm
\noindent
\textbf{Proof of Theorem \ref{1.4}.} By Lemma \ref{4.1}(ii), we have 
$$ ||\ T_f(r)=O(T_g(r)) \text{ and } ||\ T_g(r)=O(T_f(r)).$$
By changing indices if necessary, we may assume that
$$\underbrace{\dfrac{Q^{\frac{d}{d_{n+3}}}_{n+3}(f)}{Q^{\frac{d}{d_{n+3}}}_{n+3}(g)}\equiv \cdots\equiv \dfrac{Q^{\frac{d}{d_{k_1}}}_{k_1}(f)}{Q^{\frac{d}{d_{k_1}}}_{k_1}(g)}}_{\text { group } 1}\not\equiv
\underbrace{\dfrac{Q^{\frac{d}{d_{k_1+1}}}_{k_1+1}(f)}{Q^{\frac{d}{d_{k_1+1}}}_{k_1+1}(g)}\equiv \cdots\equiv\dfrac{Q^{\frac{d}{d_{k_2}}}_{k_2}(f)}{Q^{\frac{d}{d_{k_2}}}_{k_2}(g)}}_{\text { group } 2}$$
$$\not\equiv \underbrace{\dfrac{Q^{\frac{d}{d_{k_2+1}}}_{k_2+1}(f)}{Q^{\frac{d}{d_{k_2+1}}}_{k_2+1}(g)}\equiv \cdots\dfrac{Q^{\frac{d}{d_{k_3}}}_{k_3}(f)}{Q^{\frac{d}{d_{k_3}}}_{k_3}(g)}}_{\text { group } 3}\not\equiv \cdot\cdot\cdot\not\equiv \underbrace{\dfrac{Q^{\frac{d}{d_{k_{s-1}+1}}}_{k_{s-1}+1}(f)}{Q^{\frac{d}{d_{k_{s-1}+1}}}_{k_{s-1}+1}(g)}\equiv \cdots \dfrac{Q^{\frac{d}{d_{k_s}}}_{k_s}(f)}{Q^{\frac{d}{d_{k_s}}}_{k_s}(g)}}_{\text { group } s},$$
where $k_s=q.$ 

If there exist a group containing more than $[\dfrac{q-n-2}{2}]$ elements then we have the desired conclusion of the theorem. We now suppose that the number of elements of each group is at most $[\dfrac{q-n-2}{2}]$.

For each $n+3\le i \le q,$ we set
\begin{equation*}
\sigma (i)=
\begin{cases}
i+[\dfrac{q-n-2}{2}]& \text{ if $i+[\dfrac{q-n-2}{2}]\leq q$},\\
i+[\dfrac{q-n-2}{2}]-q+n+2&\text{ if  $i+[\dfrac{q-n-2}{2}]> q$},
\end{cases}
\end{equation*}
and  
$$P_i=Q^{\frac{d}{d_i}}_i(f)Q^{\frac{d}{d_{\sigma (i)}}}_{\sigma (i)}(g)-Q^{\frac{d}{d_i}}_i(g)Q^{\frac{d}{d_{\sigma (i)}}}_{\sigma (i)}(f).$$

Since the number of elements of each group is at most $[\dfrac{q-n-2}{2}]$, then $\dfrac{Q^{\frac{d}{d_i}}_i(f)}{Q^{\frac{d}{d_i}}_i(g)}$ and $\dfrac{Q^{\frac{d}{d_{\sigma (i)}}}_{\sigma (i)}(f)}{Q^{\frac{d}{d_{\sigma (i)}}}_{\sigma (i)}(g)}$ belong to two distinct groups, hence $P_i\not\equiv 0$ for every $n+3\le i\le q.$ Then we have
$$
P:=\prod_{i=n+3}^{q}P_i\not\equiv 0.
$$
We set
$$
S=\bigcup_{1\le i_1<\cdots <i_{k+1}\le n+1}\bigg (\bigcap_{j=1}^{k+1}\zero Q_{i_j}(f)\bigg ).
$$
Then $S$ is an analytic set of codimension at least $2$ of $\C^m$.

\textbf{\textit{Claim: }}\  $||\ N_{P_i}(r) \ge2\sum_{i=1}^q\dfrac{d}{d_i}N^{L_i}_{Q_i(f)}$.

Indeed, fix a point $z\not\in I(f)\cup I(g)\cup S.$ We assume that $z$ is a zero of some functions $Q_i(f)$ $(1\le i\le q)$. We set
\begin{align*}
I&=\{i \ :\ 1\le i\le n+2, (f,H_i)(z)=0\} \text{ and } t=\sharp I,\\
J&=\{i \ :\ n+3\le i\le q, (f,H_i)(z)=0\} \text{ and } l=\sharp J.
\end{align*}
Here we note that $0\le t,l\le k$ and $1\le t+l\le k$. For each index $i$, it is easy to see that
\begin{align*}
\begin{cases}
\nu_{P_i}(z)\ge \dfrac{d}{d_i}\min\{\nu^0_{Q_i(f)},L_i\}& \text{ if }i\in J,\sigma (i)\not\in J\\
\nu_{P_i}(z)\ge \dfrac{d}{d_{\sigma (i)}}\min\{\nu^0_{Q_{\sigma (i)}(f)},L_{\sigma (i)}\}& \text{ if }i\not\in J,\sigma (i)\in J\\
\nu_{P_i}(z)\ge \dfrac{d}{d_i}\min\{\nu^0_{Q_i(f)},L_i\}+\dfrac{d}{d_{\sigma (i)}}\min\{\nu^0_{Q_{\sigma (i)}(f)},L_{\sigma (i)}\}& \text{ if }i,\sigma (i)\in J\\
\nu_{P_i}(z)\ge 0&\text{ if }i,\sigma (i)\not\in J\text{ and }t=0\\
\nu_{P_i}(z)\ge 1&\text{ if }i,\sigma (i)\not\in J\text{ and }t>0.
\end{cases}
\end{align*}
We set $v(z)=\sharp \{j\ :\ j,\sigma (j)\not\in J\}.$ It easy to see that 
\begin{align*}
v(z)&\ge q-n-2-2l\ge \dfrac{t(q-n-2)}{k}.
\end{align*}
Then, we have the following two cases:\\
Case 1. $t=0$. Then
\begin{align*}
\nu_{P}(z)\ge& 2\sum_{\overset{i=n+3}{i\in J}}^q\dfrac{d}{d_i}\min\{\nu^0_{Q_i(f)},L_i\}\\
=&2\sum_{i=n+3}^q\dfrac{d}{d_i}\min\{\nu^0_{Q_i(f)},L_i\}+\dfrac{q-n-2}{k}\sum_{i=1}^{n+2}\min\{\nu^0_{Q_i(f)},1\}.
\end{align*}
Case 2. $0<t\le k.$ Then
\begin{align*}
\nu_{P}(z)\ge& 2\sum_{\overset{i=n+3}{i\in J}}^q\dfrac{d}{d_i}\min\{\nu^0_{Q_i(f)},L_i\}+v(z)\ge2\sum_{i=n+3}^q\dfrac{d}{d_i}\min\{\nu^0_{Q_i(f)},L_i\}+\dfrac{t(q-n-2)}{k}\\
=&2\sum_{i=n+3}^q\dfrac{d}{d_i}\min\{\nu^0_{Q_i(f)},L_i\}+\dfrac{q-n-2}{k}\sum_{i=1}^{n+2}\min\{\nu^0_{Q_i(f)},1\}.
\end{align*}
Therefore, from the above two cases it follows that
\begin{align*}
\nu_P(z)\ge2\sum_{i=n+3}^q\dfrac{d}{d_i}\min\{\nu^0_{Q_i(f)},L_i\}+\dfrac{q-n-2}{k}\sum_{i=1}^{n+2}\min\{\nu^0_{Q_i(f)},1\}
\end{align*}
for all $z$ outside the analytic set $I(f)\cup I(g)\cup S$. 

Integrating both sides of the above inequality, we get
\begin{align}\nonumber
N_{P}(r)\ge &2\sum_{i=n+3}^q\dfrac{d}{d_i}N^{[L_i]}_{Q_i(f)}(r)+\dfrac{q-n-2}{k}\sum_{i=1}^{n+2}N^{[1]}_{Q_i(f)}(r)\\
\label{4.6}
\ge &2\sum_{i=n+3}^q\dfrac{d}{d_i}N^{[L_i]}_{Q_i(f)}(r)+\sum_{i=1}^{n+2}\dfrac{q-n-2}{kL_i}N^{[L_i]}_{Q_i(f)}(r)\ge2\sum_{i=1}^q\dfrac{d}{d_i}N^{[L_i]}_{Q_i(f)}(r).
\end{align}
Here we note that $\dfrac{q-n-2}{kL_i}\ge \dfrac{2kL}{kL_i}\ge\dfrac{2d}{d_i}\ (1\le i\le n+2)$.

Similarly, we have 
\begin{align}\label{4.7}
N_{P}(r)\ge 2\sum_{i=1}^q\dfrac{d}{d_i}N^{[L_i]}_{Q_i(g)}(r).
\end{align}

Then by (\ref{4.6}) and (\ref{4.7}) and by Theorem A with $\epsilon =\dfrac{1}{2}$, we have
\begin{align}\label{4.8}
||\ N_{P}(r)\ge d(q-n-\dfrac{3}{2})(T_f(r)+T_g(r))+o(T_f(r)).
\end{align}
Repeating the same argument as in the proof of Theorem \ref{1.3}, by Jensen's formula and by the definition of the characteristic function, we have
\begin{align}\nonumber
||\ N_{P}(r)= \sum_{i=n+3}^qN_{P_i}(r)&\le\sum_{i=n+3}^qd(T_f(r)+T_g(r))\\
\label{4.9} 
&=d(q-n-2)(T_f(r)+T_g(r)) +o(T_f(r)).
\end{align}
From (\ref{4.8}) and (\ref{4.9}), we have
\begin{align*}
||\ d(q-n-\dfrac{3}{2})(T_f(r)+T_g(r))\le d(q-n-2)(T_f(r)+T_g(r)) +o(T_f(r)).
\end{align*}
Letting $r\longrightarrow +\infty$, we get $q-n-\dfrac{3}{2}\le q-n-2$. This is a contradiction. Therefore the supposition is impossible.

Hence there must exist a group containing more than $[\dfrac{q-n-2}{2}]$ elements, then we have the desired conclusion of the theorem. \hfill $\square$

{\it Department of Mathematics,

Hanoi National University of Education,}

{\it CAU GIAY - HANOI - VIETNAM}

E-mail: quangsd@hnue.edu.vn

\end{document}